\documentstyle[12pt]{article}
\newtheorem{theorem}{Theorem}

\newtheorem{claim}{Claim}

\newtheorem{remark}{Remark}

\newcommand{\text}[1]{\quad\mbox{#1}\quad}
\def\beq{\begin{equation}}\def\eeq{\end{equation}}
\def\beqn{\begin{eqnarray}}\def\eeqn{\end{eqnarray}}
\def\pont{\hspace{-6pt}{\bf.\ }}

\def\qed{\ifhmode\unskip\nobreak\fi\quad\ifmmode\Box\else$\Box$\fi}

\title{Replacing the host $K_n$ by $n$-chromatic graphs in Ramsey-type results}

\author{Arie Bialostocki\\[-0.8ex]
\small \texttt{arie.bialostocki@gmail.com}\and
Andr\'as Gy\'arf\'as\thanks{Research was supported in part by OTKA K104373.}\\
\small Alfr\'ed R\'enyi Institute of Mathematics\\[-0.8ex]
\small Hungarian Academy of Sciences\\[-0.8ex]
\small Budapest, P.O. Box 127 \\[-0.8ex]
\small Budapest, Hungary, H-1364 \\[-0.8ex]
\small \texttt{gyarfas@sztaki.hu}\\[-0.8ex]}
\def\boxit#1{\medskip\vbox{\hrule \hbox{\vrule\kern15pt\vbox{\kern5pt
\vbox{\advance\hsize -30pt #1\par} \kern5pt}\kern15pt\vrule}\hrule} \medskip}

\begin{document}
\maketitle
\begin{abstract}

We extend two well-known results in Ramsey theory from from $K_n$ to arbitrary $n$-chromatic graphs. The first is a note of Erd\H os and Rado  stating that in every 2-coloring of the edges of $K_n$ there is a monochromatic tree on $n$ vertices. The second is the theorem of Cockayne and Lorimer stating that  for positive integers satisfying $n_1=\max\{n_1,n_2,\dots,n_t\}$ and with $n=n_1+1+\sum_{i=1}^t (n_i-1)$, the following holds.  In every coloring of the edges of $K_n$ with colors $1,2\dots,t$ there is a monochromatic matching of size $n_i$ for some $i\in \{1,2,\dots,t\}$.
\end{abstract}

\eject
Our note stems from two research directions in combinatorics. The first one is generalization of theorems by replacing $K_n$ by $n$-chromatic graphs, the second one is Ramsey theory for sparse graphs. Among the theorems from the first direction which admit a generalization is the classical theorem of R\'edei \cite{RE} stating that every tournament on $n$ vertices contains a Hamiltonian path. This theorem extends to any $n$-chromatic digraph known as the Gallai - Roy - Vitaver theorem \cite{BG}. A recent example is a theorem of the second author \cite{GYpack} which shows that a sequence $T_1,T_2,\dots,T_{n-1}$ of trees ($T_i$ has $i$ edges) can be packed into any $n$-chromatic graph if it can be packed into $K_n$.  Among the theorems which do not admit such a generalization is the Graham-Pollack theorem \cite{GP} stating that the edge set of $K_n$ cannot be partitioned into less than $n-1$ complete bipartite graphs. The possible extension of this to any $n$-chromatic graph was suggested by Alon, Saks and Seymour but eventually this was disproved by Huang and Sudakov \cite{HS}.

The second direction is  Ramsey theory for sparse host graphs, for example where $K_n$ is replaced with graphs of bounded clique number, large girth, bounded minimum degree and $n$-chromatic graphs. We recall that the Ramsey number $R(H_1,H_2,\dots,H_t)$ is the smallest integer $n$ for which the following is true: in any coloring of the edges of the complete graph $K_n$ with colors $1,2,\dots,t$, there is a monochromatic copy of $H_i$ in color $i$ for some $i\in \{1,2,\dots,t\}$.  Furthermore, the chromatic number of a graph $G$, denoted by $\chi(G)$, is the minimum number $m$ of colors for which one can color the vertices of $G$ with $m$ colors so that no two adjacent vertices are colored with the same color.

This note explores the following question: when can we extend a Ramsey-type result by replacing the host graph $K_n$ by any $n$-chromatic graph? In particular, when does it hold that every $t$-edge-colored $R(H_1,H_2,\dots,H_t)$-chromatic graph contains a monochromatic copy of $H_i$ for some $i\in \{1,2\dots,t\}$?  A well-known result of Erd\H os states that there are graphs with arbitrary large chromatic number and arbitrary large girth, therefore positive answer to our question can be expected only if all the graphs $H_i$ are acyclic. For $t=1$, $R(H_1)=|V(H_1)|$, and the answer is indeed positive, from the folklore statement that for any acyclic $H_1$, any $|V(H_1)|$-chromatic graph contains a copy of $H_1$. We shall prove
two further results that give positive answer to our question.

A remark of Erd\H os and Rado states that in any $2$-coloring of
the edges of a complete graph $K_n$ there is a
monochromatic spanning tree. This remark have been generalized in
many ways, a survey on this subject is \cite{GYSUR}. Our extension is the following.

\begin{theorem}\label{NEW1}\pont  In any $2$-coloring of the edges of any $n$-chromatic
graph $G$ there is a monochromatic tree on $n$ vertices.
\end{theorem}

\noindent {\bf Proof.} Consider a coloring of the edges of
$G$ with two colors, say red and blue. Let $\cal{H}$ be the
hypergraph on vertex set $V(G)$ whose hyperedges are the vertex sets
of the connected components $C_i$ (in both colors). Since each
vertex of $\cal{H}$ is in one red component and in one blue
component, the dual of $\cal{H}$ is a bipartite multigraph $B$.
Observe that two vertex-disjoint edges $e,f\in E(B)$ correspond to
two vertices $v_e,v_f\in V(G)$ that are not covered by any hyperedge
(component) in $\cal{H}$, in particular $(v_e,v_f)\notin E(G)$.
Therefore  any set of pairwise disjoint edges in $B$ corresponds to
an independent set in $G$. Consequently, $\chi(G)\le \chi^{'}(B)$,
where $\chi^{'}$ is the chromatic index. By K\"onig's well-known
theorem, $\chi^{'}(B)=\Delta(B)$ where $\Delta$ is the maximum degree.  Thus
$$n=\chi(G)\le \chi^{'}(B)=\Delta(B)=\max |C_i|$$ proving the theorem. \qed

\begin{remark}\pont The proof of Theorem \ref{NEW1} would be suitable to get an extension for $t$-colorings if the following
is true (a variant of the Ryser - Lov\'asz conjecture). Any $t$-partite $t$-uniform multihypergraph $\cal{H}$ satisfies
$$\chi{'}({\cal{H}})\le (t-1)\Delta({\cal{H}}).$$
\end{remark}

\begin{remark}\pont It is natural to ask which properties of monochromatic spanning trees of $2$-colored $K_n$-s can be also ensured for the monochromatic $n$-vertex trees of $n$-chromatic graphs in Theorem \ref{NEW1}. In particular, do they have diameter 3 trees? Zero sum trees? Non-crossing trees? Brooms? (For references see \cite{GYSUR}.)
\end{remark}

Our second result is Theorem \ref{clext}, the extension of the following classical result of Cockayne and Lorimer \cite{CL}.
For positive integers satisfying $n_1=\max\{n_1,n_2,\dots,n_t\}$,
$$R(n_1K_2,n_2K_2,\dots,n_tK_2)=n_1+1+\sum_{i=1}^t (n_i-1)$$
where  $n_iK_2$ is the matching of size $n_i$, i.e. $n_i$ pairwise disjoint edges.

\begin{theorem}\pont\label{clext} Suppose that $n_1=\max\{n_1,n_2,\dots,n_t\}$ and $G$ is graph such that
$\chi(G)\ge n_1+1+\sum_{i=1}^t (n_i-1)$. Then in every coloring of the edges of $G$ with $t$ colors, there is a monochromatic $n_iK_2$
for some $i,1\le i \le t$.
\end{theorem}

\noindent {\bf Proof. } We follow the nice proof in \cite{CL}. Assume $G$ is colored with colors $1,2,\dots,t$.
Set $R=n_1+1+\sum_{i=1}^t (n_i-1)$ and consider a minimal counterexample $G$, first with respect to $t$, then with respect to the number of vertices. In particular, we have $\chi(G)=R$. We define the colored complete graph $K$ with vertex set $V(G)$ by extending the coloring of $G$ with all edges of the complement of $G$ as edges of $K$ with color $0$. Colors $1,2,\dots,t$ are called genuine colors. Let $c(uv)$ denote the color of edge $uv$ in $K$.

By assumption, $K$ has no subgraph $n_iK_2$ in genuine colors. A subgraph $H\subset K$ is {\em properly colored} if no two adjacent edges of $H$ have the same color (for colors $0,1,\dots,t$).

\begin{claim}\pont\label{cyc} There is no properly colored cycle $C\subset K$.
\end{claim}

\noindent {\bf Proof of Claim \ref{cyc}.} Suppose that $C$ is a properly colored cycle in $K$, let $m_i$ denote the number of edges of $C$ with color $i$, where $i=0,1,\dots,t$. Note that for any genuine color $i$, $m_i<n_i$. Let $G_1$ be the subgraph of $G$ induced by $V(C)$ and let $G_2$ be the subgraph of $G$ induced by $V(G)\setminus V(C)$.
Observe that $V(C)$ can be covered with $\sum_{i=1}^t m_i$ independent sets: with edges in color $0$ and with vertices.
Therefore $\chi(G_1)\le \sum_{i=1}^t m_i$. Using this,
$$\chi(G_2)\ge\chi(G)-\chi(G_1)\ge n_1+1+\sum_{i=1}^t (n_i-1)-\sum_{i=1}^t m_i=n_1+1+\sum_{i=1}^t (n_i-m_i-1)\ge$$
$$\max \{n_i-m_i: 1\le i \le t\}+1 +\sum_{i=1}^t (n_i-m_i-1).$$
Since $G_2$ cannot contain monochromatic $(n_i-m_i)P_2$ in any genuine color, $G_2$ is a counterexample, smaller than $G$, contradiction.\qed

\begin{claim}\pont\label{tree} There is no properly colored subtree  $T\subset K$ such that $T$ has at least one edge in all genuine colors.
\end{claim}
\noindent {\bf Proof of Claim \ref{tree}. } The proof is along the same line as the proof of Claim \ref{cyc}. However, replacing $C$ with a tree, we get $\chi(G_1)\le 1+\sum_{i=1}^t m_i$. Continuing the same way as before,

$$\chi(G_2)\ge\chi(G)-\chi(G_1)\ge n_1+1+\sum_{i=1}^t (n_i-1)-(1+\sum_{i=1}^t m_i)=n_1+\sum_{i=1}^t (n_i-m_i-1)\ge$$
$$\max \{n_i-m_i: 1\le i \le t\}+1 +\sum_{i=1}^t (n_i-m_i-1).$$
Note that the last inequality is still true because $m_i\ge 1$ for all $i\in \{1,2,\dots,t\}$. Again, $G_2$ is a counterexample, smaller than $G$, contradiction.\qed

Select a properly colored subtree $T\subset K$ such that $T$ has {\em as many vertices as possible.} Since any edge forms a properly colored subtree,  such a tree exists with $|V(T)|\ge 2$. We shall derive a contradiction to Claim \ref{tree} by showing that $T$ has at least one edge in all genuine colors.

Assume there is a genuine color, say color $t$, such that no edge of $T$ is colored with color $t$. However, in our counterexample some edge $uv\in E(G)$ must be colored with $t$.  From the choice of $T$ and from Claim \ref{cyc}, $\{u,v\}\cap V(T)=\emptyset$.

Let $w_1$ be an arbitrary pendant vertex of $T$. From the maximality of $T$, $c(uw_1)=c(w_1w_2)$ where $w_2$ is the neighbor of $w_1$ in $T$. Similarly, there exists $w_3\in V(T)$ such that $w_2w_3\in E(T)$ and $c(uw_2)=c(w_2w_3)$. Note that $w_3\ne w_1$ otherwise deleting $w_1$ (together with the edge $w_1w_2$) from $T$ and adding the path $w_2u,uv$ we could have a properly colored subtree larger than $T$.

Consider the longest path $P=w_1,w_2,\dots, w_j$ in $T$ such that $c(uw_i)=c(w_iw_{i+1})$ for every $i\in \{1,2\dots,j-1\}$. From the previous paragraph we know that $j\ge 3$. From the definition of $T$ and $P$ it follows that $c(uw_j)=c(w_{j-1}w_j)$. Let $S$ denote the set of neighbors of $w_{j-1}$ in $T$ and set $S_1=S\setminus \{w_{j-2},w_j\}$ ($S_1=\emptyset$ is possible). Notice that for every $s\in S_1$, we have $c(vs)=c(w_{j-1}s)$ otherwise $u,v,s,w_{j-1},u$ form a properly colored cycle.

Consider the tree $U\subset K$ obtained from $T$ as follows. First remove $w_{j-1}$ together with all edges going from $w_{j-1}$ to $S$. Then add edges  $uv,uw_{j-2},uw_j$ and the set of edges $\{vs:s\in S_1\}$. Now $|V(U)|=|V(T)|+1$ thus the coloring of $U$ cannot be proper, there exist edges $xy,xz\in E(U)$ with $c(xy)=c(xz)$.  Since the coloring of $T$ is proper, w.l.o.g. $xy\in E(U)\setminus E(T)$.

\begin{itemize}

\item 1. If $x=u$ then $y,z\in \{v,w_{j-2},w_j\}$ and $c(xy)=c(xz)$ is impossible because the colors of the edges $uv,uw_{j-2},uw_j$ are all different.

\item  2.  If $x=v$ then $y,z\in S_1\cup \{u\}$. If $(x,y)=(v,s),(x,z)=(v,s')$ for $s,s'\in S_1,s\ne s'$, then $c(xy)=c(xz)$ is impossible since $c(vs)=c(w_{j-1}s)\ne c(w_{j-1}s')=c(vs')$. If $(x,y)=(v,u),(x,z)=(v,s)$ for $s\in S_1$, then $c(vu)\ne c(vs)$ otherwise $vs$ would extend $T$ to a larger properly colored tree.

\item 3. If $(x,y)=(w_{j-2}u)$ then $xz\in E(T)$ and $c(xz)=c(xy)=c(w_{j-2}u)=c(w_{j-2}w_{j-1})$ a contradiction since $xz$ and $w_{j-2}w_{j-1}$ are edges of $T$ with the same color.

\item 4. If $(x,y)=(w_j,u)$ then $xz\in E(T)$ and with $w_{j+1}=z$, $P\cup w_{j+1}$ we could have a path longer than $P$.

\item  5. If $x\in S_1$  then $y=v$. Now $c(xy)=c(xz)$ implies that $c(xz)=c(xw_{j-1})$ which contradicts the fact that $T$ is properly colored.
\end{itemize}

We get contradiction, finishing the proof. \qed

\bigskip

\begin{remark}\pont
Zolt\'an Kir\'aly \cite{KI} noticed that Theorem \ref{clext} follows immediately using \cite{CL} as a black box. Indeed, to a partition of a $t$-colored graph $G$ into $\chi(G)$ independent sets, where $\chi(G)= n_1+1+\sum_{i=1}^t (n_i-1)=n$ one can define a $t$-colored $K_n$ by coloring edge $ij$ with any color that appears on some edge between the $i$-th and $j$-th partition classes. Applying Cockayne-Lorimer theorem, there is a matching $n_iK_2 \subset K_n$ in color $i$ and it corresponds to a monochromatic matching in $G$.
\end{remark}

\begin{remark}\pont
We are not aware of any example of acyclic graphs $H_1,\dots,H_t$ for which there exists a $R(H_1,\dots,H_t)$-chromatic $G$ with a $t$-edge-coloring that has no monochromatic $H_i$ in color $i$ for all $i\in \{1,2,\dots,t\}$.
This stimulated ``Good graph hunting'', an undergraduate research project at the Budapest Semesters of Mathematics program, supervised by the second author. An acyclic graph $H$ is {\em $t$-good} if every $t$-edge coloring  of any $R(H,H,\dots,H)$-chromatic graph contains a monochromatic copy of $H$ (there are $t$ arguments in the Ramsey function). Garrison \cite{G} proved that stars are $t$-good, as well as the path $P_4$ (except possibly for $t=3$) and that $P_5,P_6,P_7$ are $2$-good. In another research project with Riasanovsky and Sherman-Bennett \cite{GYRS}, the generalization of good graphs to (acyclic) hypergraphs was explored and many $2$-good $3$-uniform hypergraphs were found - but no bad ones at all.
\end{remark}

\end{document}